\documentclass[12pt]{article}
\usepackage{amssymb}
\usepackage{latexsym,bm}
\usepackage{graphicx}
\usepackage{amsmath}
\usepackage{mathrsfs}

\date{}
\begin{document}
\title{A conjecture implying Thomassen's chord conjecture in graph theory}
\author{\hskip -10mm Xingzhi Zhan\footnote{E-mail address: \tt zhan@math.ecnu.edu.cn}\\
{\hskip -10mm \small Department of Mathematics, East China Normal University, Shanghai 200241, China}}\maketitle
\begin{abstract}
 Thomassen's chord conjecture from 1976 states that every longest cycle in a $3$-connected graph has a chord. This is one of the most important unsolved problems
 in graph theory. We pose a new conjecture which implies Thomassen's conjecture. It involves bound vertices in a longest path between two vertices in a $k$-connected graph.
 We also give supporting evidence and analyze a special case. The purpose of making this new conjecture is to explore the surroundings of Thomassen's conjecture.
\end{abstract}

{\bf Key words.} Longest cycle; longest path; chord in a cycle; $k$-connected graph

{\bf Mathematics Subject Classification.} 05C38, 05C40, 05C35
\vskip 8mm

We consider finite simple graphs. Thomassen's famous chord conjecture from 1976 is as follows ([1, Conjecure 8.1] and [6, Conjecture 6]).

{\bf Conjecture 1.} (Thomassen, 1976) {\it Every longest cycle in a $3$-connected graph has a chord.}

Conjecture 1 has been proved to be true by Thomassen himself for cubic graphs ([7] and [8]). This beautiful and important conjecture was selected as No. 65 of the unsolved problems
in the textbook [4], as Conjecture 5 in the survey article [2] collecting problems posed by Thomassen, and in Section 13 of the survey article [3] entitled
``Beautiful conjectures in graph theory". It challenges our understanding of the structure of a graph. Since a $3$-connected graph has minimum degree at least $3,$ Conjecture 1 is implied by the following conjecture [5, p.6].

{\bf Conjecture 2.} (Harvey, 2017) {\it Every longest cycle in a $2$-connected graph with minimum degree at least $3$ has a chord.}

 We denote by $V(G)$ the vertex set of a graph $G,$ by $E(G)$ the edge set of $G,$ and by $N_G(v)$ the neighborhood in $G$ of a vertex $v.$

{\bf Definition.} Let $H$ be a subgraph of a graph $G.$ A vertex $v$ of $H$ is said to be {\it $H$-bound} if all the neighbors of $v$ in $G$ lie in $H;$ i.e., $N_G(v)\subseteq V(H).$

A vertex of a path $P$ is called an {\it internal vertex} if it is not an endpoint of $P.$ For two distinct vertices $x$ and $y,$ an {\it $(x,y)$-path} is a path whose endpoints are
$x$ and $y.$ We pose the following conjecture which implies Conjecture 2, and hence Conjecture 1.

{\bf Conjecture 3.} {\it Let $G$ be a $k$-connected graph with $k\ge 2$ and let $x,\,y$ be two distinct vertices of $G.$ If $P$ is a longest $(x,y)$-path in $G,$
then $P$ contains $k-1$ internal $P$-bound vertices.}

A computer search shows that Conjecture 3 holds for all graphs of order $\le 10,$ for cubic graphs of order $\le 18,$ for $4$-regular graphs of order $\le 14,$
for triangle-free graphs of order $\le 12$ and for $C_4$-free graphs of order $\le 13.$

The case $k=2$ of Conjecture 3 without the word ``internal" (weaker version) has the following form:

{\bf Conjecture 4.}  {\it Let $G$ be a $2$-connected graph and let $x,\,y$ be two distinct vertices of $G.$ If $P$ is a longest $(x,y)$-path in $G,$
then $P$ contains a $P$-bound vertex.}

{\bf Proof that Conjecture 4 implies Conjecture 2.} Let $C$ be a longest cycle in a $2$-connected graph $G$ with minimum degree at least $3.$ Choose two consecutive vertices $x,y$ on $C.$ Clearly $x$ and $y$ partition $C$ into two paths, the longer of which we denote by $P=C[x,\,y].$ Then $P$ is a longest $(x,y)$-path in $G.$ By Conjecture 4, $P$ contains a $P$-bound vertex $v.$ Since $v$ has degree at least $3,$ and all neighbors of $v$ lie in $C,$ it follows that $C$ has a chord that is incident to $v.$ \hfill $\Box$

Note that in the above proof we do not require that the $P$-bound vertex is an internal vertex.

The following result shows that the conclusion in Conjecture 3 holds for longest paths in the whole graph. We can even relax the connectivity condition to minimum degree.
This observation is due to Guantao Chen (Private communication in July 2023).

{\bf Theorem 5.} {\it If $Q$ is a longest path in a graph of minimum degree $d$ with $d\ge 2,$ then $Q$ contains $d-1$ internal $Q$-bound vertices.}

{\bf Proof.} Let $Q=u_1,u_2,\ldots,u_k.$ Since $Q$ is a longest path, $u_1$ is $Q$-bound. Since the minimum degree of the graph is $d,$ $u_1$ has
$d$ distinct neighbors $u_2, w_1, \ldots, w_{d-1}$ where $w_j=u_{i_j}.$ Denote $f_j=u_{i_j-1},$ the predecessor of $w_j$ on $Q.$ Clearly the $d-1$ vertices
$f_1,f_2,\ldots,f_{d-1}$ are internal vertices of $Q.$ For every $j$ with $1\le j\le d-1,$ consider the path $Q_j=Q[f_j, u_1]\cup u_1w_j\cup Q[w_j, u_k].$
Then $Q_j$ is an $(f_j, u_k)$-path with the same vertex set $V(Q)$ as $Q.$ Thus every $Q_j$ is a longest path in the graph. Consequently $f_j$ is $Q_j$-bound.
Since  $V(Q_j)=V(Q),$ $f_j$ is $Q$-bound. We have shown that $f_1, f_2, \ldots, f_{d-1}$ are $d-1$ internal $Q$-bound vertices. \hfill $\Box$

{\bf Remark.} Since the endpoints of a longest path $P$ in a graph are $P$-bound, Theorem 5 has the following corollary:
 If $P$ is a longest path in a graph of minimum degree $d$ with $d\ge 1,$ then $P$ contains $d+1$  $P$-bound vertices.

An {\it independent set} in a graph is a set of vertices no two of which are adjacent.
For a set $S$ of vertices in a graph $G,$ the {\it subgraph induced by $S,$} denoted $G[S],$ is the subgraph of $G$ whose vertex set is $S$ and whose edge set consists of all those edges of $G$ which have both endpoints in $S.$ The following conjecture is equivalent to Conjecture 4.

{\bf Conjecture 6.} (The ST conjecture) {\it Suppose that the vertex set of a graph $G$ consists of two disjoint sets $S$ and $T$ such that (1) $G[S]$ is an $(x,y)$-path $P$ and
$T$ is an independent set; (2) every vertex in $S$ has at least one neighbor in $T;$ (3) every vertex in $T$ has at least two neighbors in $S.$ Then $P$ is not a longest $(x,y)$-path in $G.$ }

{\bf Proof that Conjecture 6 is equivalent to Conjecture 4.} Suppose that Conjecture 4 holds. Let a graph $G$ and an $(x,y)$-path $P$ satisfy the three conditions in Conjecture 6.
Then $G$ is connected.

Case 1. $G$ is $2$-connected. Since $P$ has no $P$-bound vertex, by Conjecture 4 we deduce that $P$ is not a longest $(x,y)$-path in $G.$

Case 2. $G$ has connectivity $1.$ Clearly any cut-vertex of $G$ is an internal vertex of $P.$ On $P$ from $x$ to $y,$ let $z$ be the first cut-vertex of $G.$
Denote $P_1=P[x,z],$ the subpath of $P$ with endpoints $x$ and $z,$ and let $S_1=V(P_1).$ Let $T_1$ be the subset of $T$ consisting of the vertices in $T$
that have a neighbor in $S_1\setminus \{z\}.$ Observe the following facts:
(1) $T_1$ is nonempty, since $N_G(x)\cap T\subseteq T_1;$ (2) any vertex in $T_1$ has no neighbor in $S\setminus S_1;$ (3) $z$ has a neighbor in $T_1,$
since otherwise $z$ would not be the first cut-vertex. Now the subgraph $G_1=G[S_1\cup T_1]$ is $2$-connected and satisfies the three conditions in Conjecture 6 with $S$ and $T$ replaced by $S_1$ and $T_1,$ respectively. Since the $(x, z)$-path $P_1$ has no $P_1$-bound vertex, it is not a longest $(x, z)$-path in $G_1$ by Conjecture 4.
Hence, there exists an $(x, z)$-path $Q$ in $G_1$ longer than $P_1.$ Thus the $(x,y)$-path $Q\cup P[z, y]$ is longer than $P,$ implying that $P$ is not a longest $(x,y)$-path in $G.$
This proves Conjecture 6.

Conversely, suppose Conjecture 6 holds. To prove Conjecture 4, let $G$ be a $2$-connected graph, let $x,\,y$ be two distinct vertices of $G,$ and let $P$ be a longest $(x,y)$-path in $G.$ We assert that $P$ contains a $P$-bound vertex. To the contrary, assume that $P$ contains no $P$-bound vertex. Denote $S=V(P)$ and let $T$ be the set of components of the graph $G-V(P).$ We define a new graph $H$ for which $V(H)=S\cup T$ and $E(H)=E(P)\cup A$ where
$$
A=\{uC|\,u\in S,\,\,\, C\in T\,\,\, {\rm and}\,\,\,u\,\,\,{\rm has\,\,\, a\,\,\, neighbor\,\,\, in\,\,\,} C\}.
$$
Then $H[S]=P$ and $T$ is an independent set of $H.$ The assumption that $P$ contains no $P$-bound vertex implies that every vertex in $S$ has a neighbor in $T.$ Since $G$ is $2$-connected, every vertex in $T$ has at least two neighbors in $S.$ By Conjecture 6, $P$ is not a longest $(x,y)$-path in $H.$ Let $W$ be an $(x,y)$-path in $H$ that is longer
than $P.$ Necessarily $W$ contains at least one vertex in $T.$ Let $V(W)\cap T=\{C_1,\ldots,C_k\},$ and let the two neighbors of $C_i$ on $W$ be $r_i$ and $f_i,$ $i=1,\ldots,k.$
Then $r_i,f_i\in V(G),$ since $T$ is an independent set in $H.$  Viewing $C_i$ as a connected subgraph of $G,$ we see that there is an $(r_i,f_i)$-path $L_i$ in $G$ whose internal vertices lie in $C_i.$ On $W$ replacing $r_iC_if_i$ by $L_i$ for every $i=1,\ldots,k,$ we obtain an $(x,y)$-path $W^{\prime}$ in $G$ that is not shorter than $W.$ Hence $W^{\prime}$ is longer than $P,$ contradicting the condition that $P$ is a longest $(x,y)$-path in $G.$ This shows that $P$ contains a $P$-bound vertex, and thus Conjecture 4 is proved.
\hfill $\Box$

\vskip 5mm
{\bf Acknowledgement.} This research  was supported by the NSFC grant 12271170 and Science and Technology Commission of Shanghai Municipality
 grant 22DZ2229014.

\end{document}